\documentclass[11pt]{article}
\usepackage[a4paper]{geometry}

\usepackage{amsmath}
\usepackage{amssymb}
\usepackage{amsthm}
\usepackage{mathrsfs}
\usepackage{bbm}
\usepackage{empheq}
\usepackage[loose]{subfigure}
\usepackage{epsfig}
\usepackage{graphicx}
\usepackage{psfrag}
\usepackage[usenames,dvipsnames]{pstricks}
\usepackage{pst-plot}
\usepackage[colorlinks]{hyperref}
\usepackage{tabls}
\usepackage{paralist}
\usepackage{hyperref}

\DeclareSymbolFontAlphabet{\Bbb}{AMSb}

\newlength{\fixboxwidth}
\newcommand{\COMMENT}[1]{}

\newcommand{\R}{\mathbb{R}}

\newcommand{\quark}{\setbox0\hbox{$x$}\hbox to\wd0{\hss$\cdot$\hss}}

\newcommand{\overbar}[1]{\mkern 1.5mu\overline{\mkern-1.5mu#1\mkern-1.5mu}\mkern 1.5mu}

\newcommand{\e}{\varepsilon}
\newcommand{\s}{\sigma}

\newcommand{\snorm}[1] {\Vert #1 \Vert}

\newtheorem{thm}{Theorem}[section]

\newtheorem{lem}[thm]{Lemma}

\theoremstyle{definition}

\hypersetup{
    linkcolor=blue,
}

\title{Separability of reproducing kernel  spaces}
\author{Houman Owhadi and Clint Scovel
\\
California Institute of Technology
}
\date{\today}

\makeatletter
\@addtoreset{equation}{section}
\@addtoreset{figure}{section}
\@addtoreset{table}{section}
\makeatother

\renewcommand{\thefigure}{\arabic{section}.\arabic{figure}}

\makeatletter
\renewcommand{\p@subfigure}{\thefigure}
\makeatother

\newcounter{mycount}

\begin{document}

\maketitle
\abstract{
We demonstrate that a reproducing kernel Hilbert or Banach space  of functions
on a separable absolute Borel space or an analytic subset of a Polish space is separable if it
 possesses a Borel measurable feature map.
}
	


\section{Introduction}

Reproducing kernel Hilbert spaces (RKHS), see e.g.~Berlinet and Thomas-Agnan \cite{Berlinet:2004} and Steinwart and Christmann \cite[Sec.~4]{SteinwartChristmann:2008},
  are important in  Statistics and Learning Theory. Moreover, when using these spaces in  Probability and Statistics, separability has powerful effects. For example, for any separable metrizable space $X$ we have:
$\mathcal{B}(X\times X)=\mathcal{B}(X)\times \mathcal{B}(X)$ for the Borel $\s$-algebras
 \cite[Prop.~4.1.7]{Dudley:2002}, the Ky-Fan metric can be defined so as to metrize convergence
in probability \cite[Thm.~9.2.2]{Dudley:2002}, 
convergence in probability implies convergence in law \cite[Prop.~9.3.5]{Dudley:2002}, convergence in law is metrized by the Prokhorov metric \cite[Thm.~11.3.3]{Dudley:2002}, the space of probability measures with the weak topology is separable and metrizable \cite[Thm.~15.12]{AliprantisBorder:2006}, and the Kantorovich-Rubinstein and Strassen theorems have sharp forms  \cite[Thms.~11.8.2 \& 11.6.2]{Dudley:2002}. Moreover, separable Hilbert spaces are Polish so that we have all the machinery
of descriptive set theory available, regular conditional probabilities exist \cite[Thm.~10.2.2]{Dudley:2002}, Bochner integration is simple \cite[Lems.~11.37 \& 11.39]{AliprantisBorder:2006}, and all probability
measures on them are tight \cite[Thm.~69, 77-III]{DellacherieMeyer:1975}.
Most importantly,  by a classical result, see e.g.~Halmos \cite[Prob.~17]{halmos1982hilbert}, all
 separable Hilbert spaces are isomorphic with
$\ell^{2}(\mathbb{N})$.

According to Montgomery \cite{montgomery1935non}, ``Separability is a property which greatly facilitates work in metric spaces, but it may be of some interest to point out that this property has been unnecessarily assumed in the proofs of
certain theorems  concerning such spaces and concerning functions defined on them.''
 Indeed, many works do assume separability of the RKHS. For example,
Steinwart and Christmann's \cite[Thm.~7.22]{SteinwartChristmann:2008} oracle inequality for SVMs,  Christmann and Steinwart \cite[Thms.~7 \& 12]{christmann2009consistency},
\cite{christmann2010universal}, \cite{christmann2008consistency}, Steinwart and Christmann \cite{steinwart2009sparsity}, \cite{steinwart2011estimating},
De Vito, Rosasco and Toigo \cite{de2014learning},
Hable and Christmann \cite[Thm.~3.2]{hable2014estimation}, Luki{\'c}
and Beder \cite{lukic2001stochastic}, Steinwart \cite{steinwart2014convergence} and
Vovk \cite[Thm.~3]{vovk2006line}.
  De Vito,  Umanit{\`a} and Villa \cite{de2013extension} assume it in their
generalization of Mercer's theorem to matrix valued kernels,
  and Christmann, Van Messem and Steinwart \cite{christmann2009consistency} assert that Support Vector Machines (SVMs) are known to be consistent
and robust for classification and regression if they
are based on a Lipschitz continuous loss function and on a
bounded kernel with a separable reproducing kernel
Hilbert space which is dense in $L^{1}(\mu)$, where $\mu$ is the marginal distribution of the data-generating distribution.   Cambanis \cite{cambanis1975measurability} proves that a stochastic process with index set a Borel subset of a Polish space has
a measurable modification if and only if   
 the reproducing kernel corresponding to the autocorrelation function is measurable and
its corresponding RKHS is separable, and that a second order process with index set the real line is oscillatory
if and only if its RKHS is separable.
Nashed and Walter \cite{nashed1991general} require a separable RKHS in their development of
 sampling theorems for functions in reproducing kernel Hilbert spaces, and Zhang and Zhang \cite{zhang2011frames}
in reproducing kernel Banach spaces.
 Hein and Bousquet \cite{hein2004kernels}
require it and give some sufficient conditions for it.  For an example of a non-separable RKHS, see
Canu, Mary, and Rakotomamonjy \cite[Ex.~8.1.6]{Canu}.

Let us now briefly discuss the topological spaces under consideration.
A Polish space is a separable completely metrizable space and a Suslin space
a Hausdorff continuous image of a Polish space.
Following Frolik
\cite{Frolik:1970},  a metrizable  space
$X$ is said to be absolute Borel if  $X \subset Z$ is a Borel subset
for all metrizable
$Z$ for which it is a subspace. 
  Moreover, Frolik \cite{Frolik:1963} introduces {\em bianalytic spaces}
as Suslin spaces such that their complement in their \v{C}ech compactification is also Suslin and,
in Frolik \cite[Thm.~12]{Frolik:1963},
shows that a metrizable space is separable absolute Borel if and only if it is bianalytic.
On the other hand, a subset of a Polish space is called
analytic if it is Suslin.
Indeed,  the two types of spaces considered here,    
separable absolute Borel spaces and analytic subsets of  Polish spaces, are very general.
For example, for a Borel subset of a Polish space, Frolik \cite[Thm.~1]{Frolik:1970} asserts
 that it is separable absolute Borel and the famous 
Suslin theorem, see e.g.~Kechris
\cite[Thm.~13.7]{Kechris:1995} asserts that it is  analytic.  That is, they  both include
 any Borel subset of a Polish space, in particular any Borel subset of a separable Banach space, so 
any Borel, open, or closed subset of 
$\R^{n}$. Since $\R^{n}$ is Polish, it follows that this class also includes
any analytic subset of $\R^{n}$.  Counter examples to these spaces include non-separable spaces, non-metrizable spaces and
non Suslin spaces. Moreover, 
Lusin's Theorem, see e.g.~Kechris
\cite[Thm.~21.10]{Kechris:1995}, asserts that 
all analytic subsets of a Polish space are universally measurable, that is for every
$\s$-finite measure it is trapped between two Borel subsets
of the same measure. Consequently, 
 any subset of a Polish space which is not universally measurable is a counterexample.

Reproducing kernel Hilbert spaces are Hilbert spaces of  real-valued functions such that
pointwise evaluation is continuous. In their generalization to reproducing kernel Banach spaces (RKBS),
 Zhang, Xu, and Zhang
 \cite{ZhangXuZhang:2009} stipulate that a RKBS on a set $X$
 is a  reflexive Banach space of real valued functions on $X$ whose dual space is isometric with a  Banach space
of functions on $X$, such that  pointwise evaluation is continuous for both the Banach space and its dual.
They then proceed to develop the theory much along the lines of RKHSs. In particular, in
 \cite[Thm.~2]{ZhangXuZhang:2009} they show that RKBSs possess reproducing kernels. Moreover,  in
\cite[Thm.~3]{ZhangXuZhang:2009} they demonstrate that if  $\Phi:X \rightarrow \mathcal{W}$ is a map
to a  reflexive Banach space $\mathcal{W}$ and $\Phi^{*}:X \rightarrow \mathcal{W}^{*}$ is a map to its dual
such that  the linear span of the image of each map is dense, then  a RKBS is determined
with reproducing kernel  $K(x,x')=[\Phi^{*}(x),\Phi(x')]$, where $[\,,\,]$ is the dual pairing
between $ \mathcal{W}^{*}$  and $ \mathcal{W}$.  Moreover, in  \cite[Thm.~4]{ZhangXuZhang:2009} they assert
that all RKBSs possess such maps. Consequently, we refer to such maps $\Phi:X \rightarrow \mathcal{W}$ and 
 $\Phi^{*}:X \rightarrow \mathcal{W}^{*}$ as primary and secondary feature maps
for the RKBS.

 This generalization to RKBSs
has generated much interest lately, for example  see recent results of Fukumizu,  Lanckriet and Sriperumbudur \cite{fukumizu2011learning}, Zhang and Zhang \cite{zhang2012regularized},
Fasshauer, Hickernell and Ye \cite{fasshauer2015solving} in Machine Learning, in particular  
of Song, Zhang and Hickernell \cite{song2013reproducing} on sparse learning, and the recent results of Zhang and Zhang \cite{zhang2011frames}, 
 Han, Nashed, and Sun \cite{han2009sampling} and  Christensen \cite{christensen2012sampling}
 concerning sampling expansions, frames and Riesz bases in Banach spaces.

It is the purpose of this paper to establish separability for both RKHSs and RKBSs
when the domain is a separable absolute Borel space or an analytic subset of a Polish space, in particular when
it is a Borel subset of a Polish space, under the  simple assumption that the reproducing kernel space possesses a Borel measurable feature map.

\section{Main Results}
Before our main results, we review some existing results regarding the separability of RKHSs.
We will  consider both when $X$ is not a topological space and when it is. 
When $X$ is not topological,
Berlinet and Thomas-Agnan \cite[Thm.~15, pg.~33]{Berlinet:2004} shows that a RKHS  $H$ is separable if there is a countable subset
$X_{0} \subset X$ such that  $f \in H$ and $f(x)=0, x \in X_{0} $ implies that
$f=0$. Moreover,  a result of  Fortet \cite[Thm.~1.2]{Fortet:1973}
asserts that a RKHS  with kernel $k$ is separable if and only if for all $\e>0$
there exists a countable partition $B_{j}, j \in \mathbb{N}$ of $X$ such that  for all $j\in
\mathbb{N}$ and all $x_{1},x_{2} \in B_{j}$  we have
\[  k(x_{1},x_{1})+k(x_{2},x_{2})-k(x_{1},x_{2})-k(x_{2},x_{1}) <\e \, .\]

Regarding the separability of RKBSs, an 
if and only if characterization is obtained through a
 generalization of Fortet's Theorem from RKHSs to RKBSs. We  suspect the proof of our version,
  Theorem  \ref{thm_fortet},
 is  similar to  Fortet's \cite[Thm.~2.1]{Fortet:1973} for RKHSs, but  it is not written down
there.  Indeed, Fortet's result is
 a regularity condition on the pullback (pseudo) metric
\[d_{\Phi}(x_{1},x_{2}):=  \snorm{\Phi(x_{1})-\Phi(x_{2})}_{H_{1}}=\sqrt{ k(x_{1},x_{1})+k(x_{2},x_{2})-k(x
_{1},x_{2})-k(x_{2},x_{1})}
\]
 to $X$ determined by a feature map $\Phi:X \rightarrow
H_{1}$. In particular,  Fortet's condition then becomes: for all $\epsilon >0$ there exists a  countable partition
$B_{j}, j \in \mathbb{N}$ of $X$ such that 
\begin{equation}
\label{cond_fortet}
d_{\Phi}(x_{1},x_{2}) <\sqrt{\e},\quad  x_{1},x_{2} \in B_{j}, j\in  \mathbb{N} \, .
\end{equation}
We begin with a preparatory lemma  asserting that the separability of the image of the feature map implies the separability
of the corresponding RKHS or RKBS. This lemma is used in both the proof of our generalization of Fortet's
result, Theorem \ref{thm_fortet}, which is valid when $X$ is not a topological space,   and our main result
Theorem \ref{thm_main}, valid when $X$ is a separable absolute Borel space or an analytic subset of a Polish space.   
\begin{lem}
\label{lem_prep}
Consider a (RKBS) RKHS $\mathcal{K}$ of functions on a set $X$ with feature (Banach) Hilbert space $\mathcal{W}$ and (primary) feature
map $\Phi:X \rightarrow \mathcal{W}$. If $\Phi(X) \subset \mathcal{W}$ is a separable
subspace, then $\mathcal{K}$ is separable.
\end{lem}
  We can now present our generalization of Fortet's result to
 RKBSs expressed in terms of  the pseudometric space $(X,d_{\Phi})$\footnote{We would like to thank
one of the referees for pointing out the  possible connection between Fortet's condition \eqref{cond_fortet}
and a Lindel\"{o}f type condition on the pseudometric space $(X,d_{\Phi})$}.
\begin{thm}
\label{thm_fortet}
A RKBS $\mathcal{K}$ of functions on a set $X$ is separable if and only if there exists a feature
Banach space $\mathcal{W}$ and feature
map $\Phi:X \rightarrow \mathcal{W}$  such that the topological space
$(X,d_{\Phi})$  determined by the pullback pseudometric
\[d_{\Phi}(x_{1},x_{2}):=  \snorm{\Phi(x_{1})-\Phi(x_{2})}_{\mathcal{W}}, \quad x_{1}, x_{2} \in X\]
is separable.
\end{thm}

Now let us consider the case when $X$ is a topological space.
Since separability is preserved under continuous maps,
see e.g.~\cite[Thm.~16.4]{Willard}, Lemma \ref{lem_prep} implies 
  the RKBS version of Steinwart and Christmann
 \cite[Lem.~4.33]{SteinwartChristmann:2008} when combined with \cite[Lem.~4.29]{SteinwartChristmann:2008}:
A  RKBS of functions on a separable space $X$ is separable
if it has a
 continuous feature map.
Steinwart and Christmann \cite[Lem.~4.33]{SteinwartChristmann:2008} assert that
if $X$ is separable and the kernel $k$ corresponding to the RKHS $H$ is continuous,
then $H$ is separable. More generally, Steinwart and Scovel \cite[Cor.~3.6]{SteinwartScovel:2011}  show that
if there exists a finite and strictly positive Borel measure
 on $X$, then every bounded and separately continuous kernel k has a separable RKHS. 
However, to obtain our main result, our primary tool to derive separability comes from  theorems of
Stone \cite[Thm.~16, pg.~32]{Stone}, when $X$ is  separable
  absolute Borel, and
Srivastava's \cite[Thm.~4.3.8]{Srivastava}  version of
Simpson \cite{simpson1985bqo} when $X$ is an  analytic subset of  a Polish space.  It is interesting
to note that Srivastava's proof is different from Simpson's in that it does not use
Stone's Theorem \cite[Thm.~16, pg.~32]{Stone}.
\begin{lem}
\label{lem_stone}
Let $X$ be separable absolute Borel or an analytic subset of a Polish space and let $Y$ be  a metric  space, and suppose that $f:X\rightarrow Y$ is Borel measurable.
Then $f(X)\subset Y$ is separable.
\end{lem}

Steinwart and Christmann \cite[Lem.~4.25]{SteinwartChristmann:2008} shows that separate measurability of the kernel
combined with separability of the corresponding RKHS implies that the canonical feature map is measurable. Our main result is a kind of converse
 when $X$ is separable absolute Borel or an analytic subset of a Polish space.
\begin{thm}
\label{thm_main}
Let $X$ be  separable absolute Borel or an analytic subset of a Polish space
 and let $\mathcal{K}$ be a $RKHS$ with measurable feature map, or
a $RKBS$ with measurable primary feature map, of real-valued functions on $X$. Then $\mathcal{K}$ is separable.
\end{thm}

\section{Proofs}
\subsection{Proof of Lemma \ref{lem_prep}}
For RKHSs this assertion is contained in the proof of
Steinwart and Christmann
 \cite[Lem.~4.33]{SteinwartChristmann:2008}.
Roughly, the argument is that rational linear combinations are dense in the linear span of
$\Phi(X)$ and the linear span is dense
in the closed linear span in the metric defined in the proof of \cite[Thm.~4.21]{SteinwartChristmann:2008}.
For the RKBS case, since  $\Phi:X\rightarrow \mathcal{W}$ is a primary feature map it
satisfies   $\overbar{span(\Phi(X))}= \mathcal{W}$. Moreover,   
 since
 $\Phi(X)\subset \mathcal{W}$ is separable,  
 the same argument as used in the RKHS
case  shows that the closed linear span $\overbar{span(\Phi(X))}= \mathcal{W}$ is separable, so we conclude that 
 $\mathcal{W}$  is separable.
 Since $\mathcal{W}$ is reflexive 
it follows from \cite[Cor.~1.12.12]{megginson2012introduction} that $\mathcal{W}^{*}$ is separable.  Moreover, Zhang, Xu and Zhang \cite[Thm.~3]{ZhangXuZhang:2009}
implies that the dual Banach space  is
\[ \mathcal{K}^{*}:=\{[\Phi(\cdot),u^{*}]: u^{*} \in \mathcal{W}^{*}\}\]
with norm
\[\|[\Phi(\cdot),u^{*}]\|_{\mathcal{K}^*}:=\|u^{*}\|_{\mathcal{W}^{*}}\, ,\]
so that the mapping from $\mathcal{W}^{*}$ to $\mathcal{K}^{*}$ defined by
$u^{*} \mapsto [\Phi(\cdot),u^{*}]$ is an isometry. Consequently, the separability of $\mathcal{W}^{*}$
implies the separability of $\mathcal{K}^{*}$. Since $\mathcal{K}$ and therefore $\mathcal{K}^{*}$ are reflexive it follows from
\cite[Cor.~1.12.12]{megginson2012introduction} that $\mathcal{K}$ is separable.

\subsection{Proof of Theorem \ref{thm_fortet}}
Let us first demonstrate the equivalence between the separability of the RKBS and Fortet's condition 
\eqref{cond_fortet}. Then we will demonstrate the equivalence between Fortet's condition
and the separability of the  pseudometric space $(X,d_{\Phi})$.  
 Let us begin with "if". To that end, let us  show that condition \eqref{cond_fortet}               implies that
 $\Phi(X)$ is separable. Indeed, fix $\e>0$ and for each     $\frac{\e}{2^{k}}, k \in \mathbb{N}$ let $B_{j}^{k}, j \in \mathbb{N}$ denote the corresponding partition and let  $x_{j}^{k} \in B_{j}^{k}$ denote a selection.
Then the set $\Phi(x^{k}_{j}), k \in \mathcal{N}, j \in\mathcal{N}$ is countable,  and it is easy to show it is dense in $\Phi(X)$. That is, $\Phi(X)$ is separable, and the separability of $\mathcal{K}$ follows
from Lemma \ref{lem_prep}.
Now for the "only if", suppose that $\mathcal{K}$ is separable. Then  the canonical feature space $\mathcal{W}:=\mathcal{K}$ is separable,
and
  since $\mathcal{K}$ is metric, by e.g.~\cite[Thm.~16.8]{Willard} it is second countable. Therefore, since
second countability is inherited by subspaces, see e.g.~\cite[Thm.~16.2]{Willard},
  it follows for the corresponding canonical feature map $\Phi:X
\rightarrow \mathcal{K}$, that  $\Phi(X) \subset \mathcal{K}$ is second countable, and therefore,
 by e.g.~\cite[Thm.~16.9]{Willard}, it is separable.  Therefore there exists
a countable dense set $\Phi(x_{j}) \in \Phi(X), j \in \mathbb{N}$. Therefore, if for each
$\e >0$ and for each $j \in \mathbb{N}$  we define $B_{j}=\{x\in X:\snorm{\Phi(x_{j})-\Phi(x)}_{\mathcal{K}} < \frac{\e}{2}\}$, it follows that $\cup_{j \in \mathbb{N}}{B_{j}}=X$ and
$\snorm{\Phi(x_{1})-\Phi(x_{2})}_{\mathcal{K}} < \e$ for all $x_{1},x_{2}\in B_{j}$.
Therefore, we have established the equivalence between the separability of the RKBS and Fortet's condition
\eqref{cond_fortet}.

Now let us demonstrate the equivalence between Fortet's condition
and the separability of the  pseudometric space $(X,d_{\Phi})$.  
To that end,  suppose that the pseudometric space $(X,d_{\Phi})$  is separable.  Then
Willard \cite[Thm.~16.11]{Willard} asserts that in a pseudometric space, the conditions of being 
Lindel\"{o}f, second countable, and separable are equivalent. Therefore $(X,d_{\Phi})$
is Lindel\"{o}f in that every open cover has a countable subcover. For $x\in X$, let
$B_{\Phi}(x,\epsilon):=\{x'\in X:d_{\Phi}(x,x') <\epsilon\}$ denote the open ball about $x$, and
   for each $\epsilon >0$
consider the open cover $\{B_{\Phi}(x,\frac{\epsilon}{2}), x \in X\}$. Then since $(X,d_{\Phi})$
is  Lindel\"{o}f it follows there exists a countable subcover 
$\{B_{\Phi}(x,\frac{\epsilon}{2}), x \in X_{0}\}$  where $X_{0}$ is a countable. This cover satisfies
Fortet's condition \eqref{cond_fortet} for the value $\epsilon$ and since $\epsilon$  was arbitrary, it follows
that the map $\Phi:X \rightarrow \mathcal{W}$
 satisfies Fortet's condition \eqref{cond_fortet}. 
In the other direction, suppose that the map $\Phi:X \rightarrow \mathcal{W}$
 satisfies Fortet's condition \eqref{cond_fortet}.  Fix $\e>0$ and for each     $\frac{\e}{2^{k}}, k \in \mathbb{N}$ let $B_{j}^{k}, j \in \mathbb{N}$ denote the corresponding partition and let  $x_{j}^{k} \in B_{j}^{k}$ denote a selection.
Then the set $\Phi(x^{k}_{j}), k \in \mathcal{N}, j \in\mathcal{N}$ is countable,  and it is easy to show it is dense in $\Phi(X)$.  That is, for $x \in X$, 
the countable set $\{\Phi(x^{k}_{j}),  k \in \mathcal{N}, j \in\mathcal{N}\}$ comes arbitrarily close
to $ \Phi(x)$. It follows that  the countable set  $\{x^{k}_{j},  k \in \mathcal{N}, j \in\mathcal{N}\}$ 
comes arbitrarily close to $x$ in the pseudometric $d_{\Phi}$. Consequently $(X,d_{\Phi})$  is separable.

\subsection{Proof of Lemma \ref{lem_stone}}
The case when $X$ is an analytic subset of a Polish space
 follows  directly from  Srivastava \cite[Thm.~4.3.8]{Srivastava}.
 When  $X$ is separable absolute Borel,
it follows from Stone's Theorem  \cite[Thm.~16, pg.~32]{Stone}  that when $Y$ is a metric space and
 $\Phi:X \rightarrow Y$ is a  measurable bijection,
that the image $Y$  is separable.
However, when $\Phi$ is not surjective,  since $\Phi(X)\subset Y$ is a metric space, the assertion that  the metric subspace
$\Phi(X) \subset Y$ is separable follows assuming that
$\Phi$ is a measurable injection. Moreover,
 injectivity is also unnecessary. To see this, extend to the injective map
$\hat{\Phi}: X \rightarrow X \times Y$ defined by
$\hat{\Phi}(x):=\bigl(x,\Phi(x)\bigr)$. Then it follows from Hansell's  \cite[Thm.~1]{Hansell} generalization
of  Kuratowski \cite[Thm.~1, Sec.~31, VI]{kuratowski1966topology}  to the nonseparable case,
 that $\hat{\Phi}$ is measurable. To see how it is obtained, since $X$ is assumed to be separable and metrizable, it is
second countable, see e.g.~\cite[Thm.~16.11]{Willard}, so that it has
a countable base $\{G_{n}, n \in \mathbb{N}\}$ of open sets generating its topology.  Let
$W \subset X \times Y$  be open and  define
\[V_{n}=\cup{\{V: V \, \, \text{open},\, G_{n}\times V \subset W\} }\, .\]
Then \[W=\mathop{\cup}_{n\in \mathbb{N}}{G_{n}\times V_{n}}\]
 and therefore
\[\hat{\Phi}^{-1}(W)=\mathop{\cup}_{n\in \mathbb{N}}{G_{n}\cap \Phi^{-1}(V_{n})}\, .\]
Since $G_{n}$ and $V_{n}$ are open and therefore measurable and $\Phi$ is measurable it follows
that $\hat{\Phi}^{-1}(W)$ is measurable. Consequently, since the open sets generate the
Borel $\s$-algebra, it follows that $\hat{\Phi}$ is Borel measurable.
Moreover, since $\hat{\Phi}$ is injective the above discussion shows that
$\hat{\Phi}(X)\subset X \times Y$ is separable. Since
$\Phi(X)= P_{Y}\hat{\Phi}(X)$ where $P_{Y}$ is the projection
onto the second component and $P_{Y}$ is continuous, and separability is preserved under continuous maps,
see e.g.~\cite[Thm.~16.4]{Willard},
it follows that $\Phi(X)\subset Y$ is separable.

\subsection{Proof of Theorem \ref{thm_main}}
Since the feature space is metric,
  Lemma \ref{lem_stone} implies that the image $\Phi(X)$ is separable for any measurable feature map
$\Phi$. The
 assertion  then follows from Lemma \ref{lem_prep}.

\section*{Acknowledgments}
The authors would like to thank the reviewers for comments and suggestions which substantially improved
both the content and the presentation of this work.
The authors gratefully acknowledge this work supported by the Air Force Office of Scientific Research under Award Number
FA9550-12-1-0389 (Scientific Computation of Optimal Statistical Estimators).

\addcontentsline{toc}{section}{References}
\bibliographystyle{plain}
\bibliography{refs}

\end{document}